# On the construction of numerical iterative schemes of any order of convergence for solving nonlinear systems of equations

Stefan HOTHAZIE*[,1], Camelia MUNTEANU[1], Mihaela NASTASE[1]

*Corresponding author
[1]INCAS – National Institute for Aerospace Research "Elie Carafoli",
B-dul Iuliu Maniu 220, Bucharest 061126, Romania,
hothazie.stefan@incas.ro*, munteanu.camelia@incas.ro, nastase.mihaela@incas.ro

***Abstract:*** *This paper presents a methodology for constructing iterative schemes of any order of convergence for solving nonlinear systems of equations. It also provides formulas for the order of convergence of any iterative schemes constructed using the method proposed in this paper. A test case is conducted numerically for the second and third order of convergence using a computer algebra system called Maxima. The code used is listed at the end of the test case.*

***Key Words:*** *nonlinear algebraic equations, iterative schemes, order of convergence, numerical applications*

## 1. INTRODUCTION

The problem of solving a nonlinear system of equations has applications in many fields of study. In the field of structural mechanics, the fundamental conservation laws are based on partial differential equations that are usually nonlinear. When converted into numerical schemes for solving complex problems, they retain their nonlinearity. In order to make this problem tractable, either of the following two techniques are employed: linearizing the numerical scheme for finding a solution using linear algebra; and using an iterative numerical scheme that directly solves the nonlinear system of equations.

There are three important practical properties that apply to numerical methods used in the latter technique.

The first property, the order of convergence, is of importance because it represents how fast the number of correct digits of the current solution grows. For example, the most commonly used iterative method is the Newton-Raphson method. It has a convergence of order two. That means that the solution's number of correct digits doubles each iteration. Another iterative method is Halley's method. It has a convergence of order three, which means that the number of exact digits triples with each iteration, and so on.

The second property is the iteration computational cost, or more precisely, how much computer time does it take to compute one iteration. The advantage of having a higher convergence rate is that the number of exact digits grows much faster. The disadvantage is that it takes a higher amount of computer time to compute each iteration. This problem can be solved in various ways, from rewriting the numerical scheme in discrete form to using mathematical manipulation to reduce the amount of calculations needed.

Lastly, the third property is the dimensionality of the numerical method: whether a single equation of one variable is solved, or a system of equations of multiple variables is solved. This property simply distinguishes the iterative methods by their ability to solve for one variable or for multiple variables.

Whenever an iterative numerical method is created, two issues must be clarified. The iteration formula, convergence rate and computational cost must be mentioned, and, if possible, a proof of the convergence rate must be given. So far, many iterative methods have been created, many of which have a proof for convergence. Unfortunately, to the best of our knowledge, there is no unifying theory underlying the domain of numerical iterative schemes for solving nonlinear systems of equations. This paper will attempt to provide such a theory with which iterative schemes of any order of convergence can be constructed. It will also present formulas for calculating order of convergence of these numerical methods. More details on the state of the art on this subject are presented in chapter 7.

## 2. ITERATIVE METHODS CONCERNING THE ONE-DIMENSIONAL CASE

The technique presented in this paper, for creating iterative schemes, focuses mainly on systems of nonlinear equations, but, for the purpose of continuity, we shall start by discussing the one-dimensional case. This reduces to the classical single variable nonlinear equation of the form:

$$f(x) = 0$$

The main idea behind obtaining the solution $\bar{x}$ is to find the inverse $f^{-1}(x)$ and evaluate it at the point $x = 0$:

$$\bar{x} = f^{-1}(0)$$

For the purpose of readability, the inverse function will be renamed as:

$$g(x) := f^{-1}(x)$$

$$g\big(f(x)\big) = f^{-1}\big(f(x)\big) = x \tag{1}$$

The steps to assembling the iterative formulas are:

1. Write the Taylor series of the inverse function in $f(x)$ and,
2. Transform all the derivatives into functions of $x$.

We start by expanding the function $g$ into its Taylor series in $f(x)$, around a point $f(x_0)$:

$$g(f(x)) = g\big(f(x_0)\big) + \frac{dg}{df}\bigg|_{f(x_0)} \big(f(x) - f(x_0)\big) + \frac{d^2g}{df^2}\bigg|_{f(x_0)} \frac{(f(x) - f(x_0))^2}{2} + \cdots$$

The next step is to evaluate the derivatives of $g$ as functions of $x$. This is done by differentiating with respect to $x$ equation (1). The following calculations are for the first derivative of $g$:

$$\frac{d}{dx}\Big(g\big(f(x)\big)\Big) = \frac{df}{dx}\bigg|_x \frac{d}{df}\Big(g\big(f(x)\big)\Big) = \frac{df}{dx}\bigg|_x \frac{dg}{df}\bigg|_{f(x)}$$

$$\frac{d}{dx}\Big(g\big(f(x)\big)\Big) = \frac{d}{dx}(x) = 1$$

$$\Rightarrow \left.\frac{dg}{df}\right|_{f(x)} = \frac{1}{\left.\frac{df}{dx}\right|_{x}} \tag{2}$$

For derivatives of higher order, we repeat the same procedure, differentiating each time the previous expression. We will do one more for the second derivative, as an example:

$$\frac{d}{dx}\left(\left.\frac{dg}{df}\right|_{f(x)}\right) = \left.\frac{df}{dx}\right|_{x}\frac{d}{df}\left(\left.\frac{dg}{df}\right|_{f(x)}\right) = \left.\frac{df}{dx}\right|_{x}\left.\frac{d^2g}{df^2}\right|_{f(x)}$$

$$\frac{d}{dx}\left(\left.\frac{dg}{df}\right|_{f(x)}\right) = \frac{d}{dx}\left(\frac{1}{\left.\frac{df}{dx}\right|_{x}}\right) = -\frac{\left.\frac{d^2f}{dx^2}\right|_{x}}{\left(\left.\frac{df}{dx}\right|_{x}\right)^2}$$

$$\Rightarrow \left.\frac{d^2g}{df^2}\right|_{f(x)} = -\frac{\left.\frac{d^2f}{dx^2}\right|_{x}}{\left(\left.\frac{df}{dx}\right|_{x}\right)^3} \tag{3}$$

We are now ready for the complete substitution. The resulting Taylor series using equations (1), (2) and (3) is:

$$x = x_0 + \frac{1}{\left.\frac{df}{dx}\right|_{x_0}}\big(f(x) - f(x_0)\big) - \frac{\left.\frac{d^2f}{dx^2}\right|_{x_0}}{\left(\left.\frac{df}{dx}\right|_{x_0}\right)^3}\frac{\big(f(x) - f(x_0)\big)^2}{2} + \cdots \tag{4}$$

Writing (4) as a fixed-point iteration, around the solution $\bar{x}$, and noting that $f(\bar{x}) = 0$ , we have:

$$\bar{x} = x_0 + \frac{1}{\left.\frac{df}{dx}\right|_{x_0}}\big(f(\bar{x}) - f(x_0)\big) - \frac{\left.\frac{d^2f}{dx^2}\right|_{x_0}}{\left(\left.\frac{df}{dx}\right|_{x_0}\right)^3}\frac{\big(f(\bar{x}) - f(x_0)\big)^2}{2} + \cdots \tag{5}$$

$$x_{n+1} = x_n + \frac{1}{\left.\frac{df}{dx}\right|_{x_n}}\big(-f(x_n)\big) - \frac{\left.\frac{d^2f}{dx^2}\right|_{x_n}}{\left(\left.\frac{df}{dx}\right|_{x_n}\right)^3}\frac{\big(-f(x_n)\big)^2}{2} + \cdots \tag{6}$$

If we had the complete infinite series, we would have got the exact solution $\bar{x}$ directly by evaluating expression (5) just once. Due to the fact that we cannot, in practice, calculate all the terms of the Taylor series expansion of the inverse function, we will be forced to truncate it after a finite number of terms and use it iteratively as a fixed-point method, like (6), to converge to a solution.

If we truncate the series after the first term we obtain the second order Newton-Raphson method:

$$x_{n+1} = x_n - \frac{f(x_n)}{\left.\frac{df}{dx}\right|_{x_n}}$$

Furthermore, if we truncate it after the third term we obtain a variant of Halley's method [4], which is a third order method:

$$x_{n+1} = x_n - \frac{f(x_n)}{\left.\frac{df}{dx}\right|_{x_n}} - \frac{f(x_n)^2 \left.\frac{d^2 f}{dx^2}\right|_{x_n}}{2\left(\left.\frac{df}{dx}\right|_{x_n}\right)^3}$$

If we continue to truncate the series after more and more terms, iteration formulas of higher and higher order will be obtained. This problem has already been solved extensively in the one-dimensional case by Schröder [2], who arrives at the same generalized iteration schemes, by using a different derivation.

## 3. CALCULATING THE ORDER OF CONVERGENCE IN THE ONE-DIMENSIONAL CASE

Let $\bar{x}$ be a solution to $f(x) = 0$ and $x_n$ a point close to $\bar{x}$. Let us rewrite the previous Taylor series (5) in compact form:

$$\bar{x} = x_n + \sum_{i=1}^{\infty} a_i \left(f(\bar{x}) - f(x_n)\right)^i \tag{7}$$

$$a_i = \frac{\left.\frac{d^i g}{df^i}\right|_{f(x_n)}}{i!}$$

The terms $a_i$ are the derivatives of the inverse function of the i$^{th}$ order, which can all be evaluated to dependent only on the variable $x$ using the procedure described in the previous chapter.

Since we deal with a truncated version of the above series, we will write for $x_{n+1}$ the first $k$ terms of the summation:

$$x_{n+1} = x_n + \sum_{i=1}^{k} a_i \left(f(\bar{x}) - f(x_n)\right)^i \tag{8}$$

If we subtract (8) from (7), we have:

$$\bar{x} - x_{n+1} = a_{k+1}\left(f(\bar{x}) - f(x_n)\right)^{k+1} + \cdots$$

We divide and multiply by $(\bar{x} - x_n)^i$ to obtain:

$$\bar{x} - x_{n+1} = a_{k+1} \left( \frac{f(\bar{x}) - f(x_n)}{\bar{x} - x_n} \right)^{k+1} (\bar{x} - x_n)^{k+1} + \cdots \tag{9}$$

Since the point $x_n$ is close to $\bar{x}$, the term in parenthesis can be further simplified:

$$\frac{f(\bar{x}) - f(x_n)}{\bar{x} - x_n} \sim \frac{f(\bar{x}) - \left( f(\bar{x}) + \left.\frac{df}{dx}\right|_{\bar{x}} (x_n - \bar{x}) + O((x_n - \bar{x})^2) \right)}{\bar{x} - x_n} = \left.\frac{df}{dx}\right|_{\bar{x}} + O\big((\bar{x} - x_n)\big) \tag{10}$$

We define the error $\delta$ as being the difference between the solution $\bar{x}$ and the current point of evaluation $x$:

$$\delta_n = \bar{x} - x_n$$

$$\delta_{n+1} = \bar{x} - x_{n+1}$$

Substituting the expressions above in (9) and (10), we have:

$$\delta_{n+1} = a_{k+1} \left( \left.\frac{df}{dx}\right|_{\bar{x}} + O(\delta_n) \right)^{k+1} \delta_n^{k+1} + \cdots = a_{k+1} \left( \left.\frac{df}{dx}\right|_{\bar{x}} \right)^{k+1} \delta_n^{k+1} + O\big(\delta_n^{k+2}\big)$$

After dropping the higher order terms, it can be clearly seen from (11),/ that a scheme, that uses $k$ terms, has an order of convergence of $k+1$:

$$\delta_{n+1} = a_{k+1} \left( \frac{df}{dx}(\bar{x}) \right)^{k+1} \delta_n^{k+1} \tag{11}$$

The conditions for $k+1$ order of convergence are the following:

- $\left.\frac{df}{dx}\right|_{x_i} \neq 0$ for all $x_i$, $i = \overline{1,n}$
- $a_{k+1}|_{x_i}$ is continuous for all $x_i$, $i = \overline{1,n}$
- $x_i$ is sufficiently close to the root $\bar{x}$

The first condition allows for the construction of the iterative formulas, because we need to divide by the first derivative to obtain the terms $a_k$, see (2) and (3).

The second condition ensures that the terms $a_j\big|_{x_i}$ for $j = \overline{1,k}$ are computable and bounded for each iteration at every point $x_i$ for $i = \overline{1,n}$. The last condition ensures $k+1$ order of convergence, because, for this to be possible, the Taylor series approximation (11) has to be valid, which involves higher order terms to be negligible.

## 4. ITERATIVE METHODS CONCERNING THE MULTIDIMENSIONAL CASE

In order to construct iteration schemes in the case of systems of nonlinear equations we have to define the notion of the inverse of a multivariate function. Since the purpose of an inverse function, when applied to the function itself, is to output the independent variable, the logical extension to multivariate functions is to have as many inverse functions as there are independent variables. Before providing the technique for the general case, we shall first start with a smaller example of a system of nonlinear equations in two variables $x_0$ and $x_1$:

$$\begin{aligned} f_0(x_0, x_1) &= 0 \\ f_1(x_0, x_1) &= 0 \end{aligned} \tag{12}$$

Since we have two variables, we now need two inverse functions $g_0$ and $g_1$, with the following properties:

$$g_0\big(f_0(x_0, x_1), f_1(x_0, x_1)\big) = x_0$$
$$g_1\big(f_0(x_0, x_1), f_1(x_0, x_1)\big) = x_1$$

We follow the same procedure as before and write the Taylor series of the two inverse functions $g_0$ and $g_1$ around a point $(x_0, x_1)$:

$$\begin{aligned} &g_0\big(f_0(\overline{x_0}, \overline{x_1}), f_1(\overline{x_0}, \overline{x_1})\big) \\ &= g_0\big(f_0(x_0, x_1), f_1(x_0, x_1)\big) + \frac{\partial g_0}{\partial f_0}\bigg|_{f_0(x_0,x_1), f_1(x_0,x_1)} \big(f_0(\overline{x_0}, \overline{x_1}) - f_0(x_0, x_1)\big) \\ &+ \frac{\partial g_0}{\partial f_1}\bigg|_{f_0(x_0,x_1), f_1(x_0,x_1)} \big(f_1(\overline{x_0}, \overline{x_1}) - f_1(x_0, x_1)\big) \\ &+ \frac{\partial^2 g_0}{\partial f_0^2}\bigg|_{f_0(x_0,x_1), f_1(x_0,x_1)} \frac{\big(f_0(\overline{x_0}, \overline{x_1}) - f_0(x_0, x_1)\big)^2}{2} \\ &+ \frac{\partial^2 g_0}{\partial f_0 \partial f_1}\bigg|_{f_0(x_0,x_1), f_1(x_0,x_1)} \big(f_0(\overline{x_0}, \overline{x_1}) - f_0(x_0, x_1)\big)\big(f_1(\overline{x_0}, \overline{x_1}) - f_1(x_0, x_1)\big) \\ &+ \frac{\partial^2 g_0}{\partial f_1^2}\bigg|_{f_0(x_0,x_1), f_1(x_0,x_1)} \frac{\big(f_1(\overline{x_0}, \overline{x_1}) - f_1(x_0, x_1)\big)^2}{2} + \cdots \end{aligned}$$

$$\begin{aligned} &g_1\big(f_0(\overline{x_0}, \overline{x_1}), f_1(\overline{x_0}, \overline{x_1})\big) \\ &= g_1\big(f_0(x_0, x_1), f_1(x_0, x_1)\big) + \frac{\partial g_1}{\partial f_0}\bigg|_{f_0(x_0,x_1), f_1(x_0,x_1)} \big(f_0(\overline{x_0}, \overline{x_1}) - f_0(x_0, x_1)\big) \\ &+ \frac{\partial g_1}{\partial f_1}\bigg|_{f_0(x_0,x_1), f_1(x_0,x_1)} \big(f_1(\overline{x_0}, \overline{x_1}) - f_1(x_0, x_1)\big) \\ &+ \frac{\partial^2 g_1}{\partial f_0^2}\bigg|_{f_0(x_0,x_1), f_1(x_0,x_1)} \frac{\big(f_0(\overline{x_0}, \overline{x_1}) - f_0(x_0, x_1)\big)^2}{2} \\ &+ \frac{\partial^2 g_1}{\partial f_0 \partial f_1}\bigg|_{f_0(x_0,x_1), f_1(x_0,x_1)} \big(f_0(\overline{x_0}, \overline{x_1}) - f_0(x_0, x_1)\big)\big(f_1(\overline{x_0}, \overline{x_1}) - f_1(x_0, x_1)\big) \\ &+ \frac{\partial^2 g_1}{\partial f_1^2}\bigg|_{f_0(x_0,x_1), f_1(x_0,x_1)} \frac{\big(f_1(\overline{x_0}, \overline{x_1}) - f_1(x_0, x_1)\big)^2}{2} + \cdots \end{aligned}$$

The next step is the computation of the derivatives of $g_0$ and $g_1$ as functions of $x_0$ and $x_1$. As an example, we will carry out the computations only for the first and second derivatives:

$$\frac{\partial}{\partial x_0}\Big(g_0\big(f_0(x_0, x_1), f_1(x_0, x_1)\big)\Big) = \frac{\partial f_0}{\partial x_0}\frac{\partial g_0}{\partial f_0}\bigg|_{f_0(x_0,x_1), f_1(x_0,x_1)} + \frac{\partial f_1}{\partial x_0}\frac{\partial g_0}{\partial f_1}\bigg|_{f_0(x_0,x_1), f_1(x_0,x_1)}$$
$$\frac{\partial}{\partial x_1}\Big(g_0\big(f_0(x_0, x_1), f_1(x_0, x_1)\big)\Big) = \frac{\partial f_0}{\partial x_1}\frac{\partial g_0}{\partial f_0}\bigg|_{f_0(x_0,x_1), f_1(x_0,x_1)} + \frac{\partial f_1}{\partial x_1}\frac{\partial g_0}{\partial f_1}\bigg|_{f_0(x_0,x_1), f_1(x_0,x_1)}$$

$$\frac{\partial}{\partial x_0}\Big(g_1\big(f_0(x_0,x_1),f_1(x_0,x_1)\big)\Big)=\frac{\partial f_0}{\partial x_0}\frac{\partial g_1}{\partial f_0}\bigg|_{f_0(x_0,x_1),f_1(x_0,x_1)}+\frac{\partial f_1}{\partial x_0}\frac{\partial g_1}{\partial f_1}\bigg|_{f_0(x_0,x_1),f_1(x_0,x_1)}$$

$$\frac{\partial}{\partial x_1}\Big(g_1\big(f_0(x_0,x_1),f_1(x_0,x_1)\big)\Big)=\frac{\partial f_0}{\partial x_1}\frac{\partial g_1}{\partial f_0}\bigg|_{f_0(x_0,x_1),f_1(x_0,x_1)}+\frac{\partial f_1}{\partial x_1}\frac{\partial g_1}{\partial f_1}\bigg|_{f_0(x_0,x_1),f_1(x_0,x_1)}$$

$$\frac{\partial}{\partial x_0}\Big(g_0\big(f_0(x_0,x_1),f_1(x_0,x_1)\big)\Big)=\frac{\partial}{\partial x_0}(x_0)=1$$
$$\frac{\partial}{\partial x_1}\Big(g_0\big(f_0(x_0,x_1),f_1(x_0,x_1)\big)\Big)=\frac{\partial}{\partial x_1}(x_0)=0$$
$$\frac{\partial}{\partial x_0}\Big(g_1\big(f_0(x_0,x_1),f_1(x_0,x_1)\big)\Big)=\frac{\partial}{\partial x_0}(x_1)=0$$
$$\frac{\partial}{\partial x_1}\Big(g_1\big(f_0(x_0,x_1),f_1(x_0,x_1)\big)\Big)=\frac{\partial}{\partial x_1}(x_1)=1$$

Written as a matrix equation, we have:

$$\begin{pmatrix}\frac{\partial g_0}{\partial f_0} & \frac{\partial g_0}{\partial f_1}\\ \frac{\partial g_1}{\partial f_0} & \frac{\partial g_1}{\partial f_1}\end{pmatrix}\begin{pmatrix}\frac{\partial f_0}{\partial x_0} & \frac{\partial f_0}{\partial x_1}\\ \frac{\partial f_1}{\partial x_0} & \frac{\partial f_1}{\partial x_1}\end{pmatrix}=\begin{pmatrix}1 & 0\\ 0 & 1\end{pmatrix}$$

Therefore, the resulting derivatives are:

$$\begin{pmatrix}\frac{\partial g_0}{\partial f_0} & \frac{\partial g_0}{\partial f_1}\\ \frac{\partial g_1}{\partial f_0} & \frac{\partial g_1}{\partial f_1}\end{pmatrix}=J^{-1}\text{ ,where }J=\begin{pmatrix}\frac{\partial f_0}{\partial x_0} & \frac{\partial f_0}{\partial x_1}\\ \frac{\partial f_1}{\partial x_0} & \frac{\partial f_1}{\partial x_1}\end{pmatrix}$$

The matrix $J$ is the Jacobian of the system. In order to compute the second derivatives, we must employ a compact system of notation that uses the Einstein summation convention. We shall use the symbol $*$, which signifies that the whole range of possible variables are taken into account. It basically shortens the following expressions:

$$x_* := x_0, x_1, x_2, x_3, x_4, \dots$$
$$f_* := f_0, f_1, f_2, f_3, f_4, \dots$$

The nonlinear system of equations (12), rewritten in the new system of notation, is:

$$f_i(x_*)=0$$

Furthermore, the previous computation of the first derivative becomes:

$$g_i\big(f_*(x_*)\big)=x_i$$

$$\frac{\partial}{\partial x_j}\Big(g_i\big(f_*(x_*)\big)\Big)=\frac{\partial f_k}{\partial x_j}\frac{\partial}{\partial f_k}\Big(g_i\big(f_*(x_*)\big)\Big)=\frac{\partial f_k}{\partial x_j}\frac{\partial g_i}{\partial f_k}\bigg|_{f_*(x_*)}$$

$$\frac{\partial}{\partial x_j}\Big(g_i\big(f_*(x_*)\big)\Big)=\frac{\partial}{\partial x_j}(x_i)=\delta_{ij}$$

$$\left.\frac{\partial f_k}{\partial x_j}\frac{\partial g_i}{\partial f_k}\right|_{f_*(x_*)} = \delta_{ij}$$

$$\left.\frac{\partial g_i}{\partial f_k}\right|_{f_*(x_*)} = \left(\frac{\partial f_k}{\partial x_j}\right)^{-1}\delta_{ij} = \left(\frac{\partial f_k}{\partial x_i}\right)^{-1} = J_{ki}^{-1}$$

The term $\delta_{ij}$ is the Kroneker delta symbol and $J_{ki}^{-1}$ is the inverse of the Jacobian matrix. We can now compute the second derivative:

$$\frac{\partial}{\partial x_p}\left(\left.\frac{\partial g_i}{\partial f_k}\right|_{f_*(x_*)}\right) = \frac{\partial f_r}{\partial x_p}\frac{\partial}{\partial f_r}\left(\left.\frac{\partial g_i}{\partial f_k}\right|_{f_*(x_*)}\right) = \left.\frac{\partial f_r}{\partial x_p}\frac{\partial^2 g_i}{\partial f_k \partial f_r}\right|_{f_*(x_*)}$$

$$\frac{\partial}{\partial x_p}\left(\left.\frac{\partial g_i}{\partial f_k}\right|_{f_*(x_*)}\right) = \frac{\partial}{\partial x_p}\left(J_{ki}^{-1}\right)$$

$$\left.\frac{\partial f_r}{\partial x_p}\frac{\partial^2 g_i}{\partial f_k \partial f_r}\right|_{f_*(x_*)} = \frac{\partial}{\partial x_p}\left(J_{ki}^{-1}\right)$$

$$\left.\frac{\partial^2 g_i}{\partial f_k \partial f_r}\right|_{f_*(x_*)} = \left(\frac{\partial f_r}{\partial x_p}\right)^{-1}\frac{\partial}{\partial x_p}\left(J_{ki}^{-1}\right) = J_{rp}^{-1}\frac{\partial J_{ki}^{-1}}{\partial x_p}$$

The Taylor series around the point $(x_*^n)$ evaluated at the solution $(\overline{x_*})$, before substitution is:

$$g_i(f_*(\overline{x_*})) = g_i(f_*(x_*^n)) + \left.\frac{\partial g_i}{\partial f_j}\right|_{f_*(x_*)}\left(f_j(\overline{x_*}) - f_j(x_*^n)\right)$$
$$+\frac{1}{2}\left.\frac{\partial^2 g_i}{\partial f_j \partial f_k}\right|_{f_*(x_*)}\left(f_j(\overline{x_*}) - f_j(x_*^n)\right)\left(f_k(\overline{x_*}) - f_k(x_*^n)\right) + \cdots$$

After substitution, and change of dummy indices, the previous Taylor series becomes the following iteration scheme:

$$x_i^{n+1} = x_i^n + J_{ji}^{-1}\left(-f_j(x_*^n)\right) + \frac{1}{2}J_{kp}^{-1}\frac{\partial J_{ji}^{-1}}{\partial x_p}\left(-f_j(x_*^n)\right)\left(-f_k(x_*^n)\right) + \cdots \quad (13)$$

If we truncate (13) to just the first two terms, we obtain the second order Newton-Raphson scheme for systems of nonlinear equations:

$$x_i^{n+1} = x_i^n - J_{ji}^{-1} f_j(x_*^n) \quad (14)$$

Truncating to the first three terms, we obtain the equivalent of Schröder's third order method for multivariate systems:

$$x_i^{n+1} = x_i^n - J_{ji}^{-1} f_j(x_*^n) + \frac{1}{2}J_{kp}^{-1}\frac{\partial J_{ji}^{-1}}{\partial x_p} f_j(x_*^n) f_k(x_*^n) \quad (15)$$

Further truncating the series to more terms, we obtain even higher orders of convergence.

# 5. CALCULATING THE ORDER OF CONVERGENCE CONCERNING THE MULTIDIMENSIONAL CASE

Let the Taylor series be rewritten in compact form as:

$$\begin{aligned}\overline{x_i} = x_i^n &+ a_{ij_0}\left(f_{j_0}(\overline{x_*}) - f_{j_0}(x_*^n)\right) + a_{ij_0j_1}\left(f_{j_0}(\overline{x_*}) - f_{j_0}(x_*^n)\right)\left(f_{j_1}(\overline{x_*}) - f_{j_1}(x_*^n)\right)\\ &+ a_{ij_0j_1j_2}\left(f_{j_0}(\overline{x_*}) - f_{j_0}(x_*^n)\right)\left(f_{j_1}(\overline{x_*}) - f_{j_1}(x_*^n)\right)\left(f_{j_2}(\overline{x_*}) - f_{j_2}(x_*^n)\right)\\ &+ \cdots\end{aligned}$$

$$\overline{x_i} = x_i^n + \sum_{p=0}^{\infty} a_{ij_0j_1\dots j_p} \prod_{k=0}^{p}\left(f_{j_k}(\overline{x_*}) - f_{j_k}(x_*^n)\right) \tag{16}$$

$$x_i^{n+1} = x_i^n + \sum_{p=0}^{m} a_{ij_0j_1\dots j_p} \prod_{k=0}^{p}\left(f_{j_k}(\overline{x_*}) - f_{j_k}(x_*^n)\right) \tag{17}$$

We define the error $\delta_*$ as being the difference between the solution $\overline{x_*}$ and the current point of evaluation $x_*$ as:

$$\begin{aligned}\delta_i^n &= \overline{x_i} - x_i^n\\ \delta_i^{n+1} &= \overline{x_i} - x_i^{n+1}\end{aligned} \tag{18}$$

Let $x_*^n$ and $x_*^{n+1}$ be two points close to the solution $\overline{x_*}$. Multiplying and dividing the term in parenthesis in (16) and (17) by $\overline{x_i} - x_i^n$, we have:

$$\overline{x_i} = x_i^n + \sum_{p=0}^{\infty} a_{ij_0j_1\dots j_p} \prod_{k=0}^{p}\left(\frac{f_{j_k}(\overline{x_*}) - f_{j_k}(x_*^n)}{\overline{x_{j_k}} - x_{j_k}^n}\right)\left(\overline{x_{j_k}} - x_{j_k}^n\right) \tag{19}$$

$$x_i^{n+1} = x_i^n + \sum_{p=0}^{m} a_{ij_0j_1\dots j_p} \prod_{k=0}^{p}\left(\frac{f_{j_k}(\overline{x_*}) - f_{j_k}(x_*^n)}{\overline{x_{j_k}} - x_{j_k}^n}\right)\left(\overline{x_{j_k}} - x_{j_k}^n\right) \tag{20}$$

The term in parenthesis from (19) and (20) can be rewritten as:

$$\begin{aligned}\frac{f_{j_k}(\overline{x_*}) - f_{j_k}(x_*^n)}{\overline{x_{j_k}} - x_{j_k}^n} &\sim \frac{f_{j_k}(\overline{x_*}) - \left(f_{j_k}(\overline{x_*}) + \left.\frac{\partial f_{j_k}}{\partial x_{j_k}}\right|_{f_*(x_*)}\left(x_{j_k}^n - \overline{x_{j_k}}\right) + O\left(\left(x_{j_k}^n - \overline{x_{j_k}}\right)^2\right)\right)}{\overline{x_{j_k}} - x_{j_k}^n}\\ &= \left.\frac{\partial f_{j_k}}{\partial x_{j_k}}\right|_{f_*(x_*)} + O\left(\left(\overline{x_{j_k}} - x_{j_k}^n\right)\right)\end{aligned} \tag{21}$$

Substituting the expression for the term in parenthesis (21) in (19) and (20), we have:

$$\overline{x_i} = x_i^n + \sum_{p=0}^{\infty} a_{ij_0j_1\dots j_p} \prod_{k=0}^{p}\left.\frac{\partial f_{j_k}}{\partial x_{j_k}}\right|_{f_*(x_*)}\left(\overline{x_{j_k}} - x_{j_k}^n\right)$$

$$x_i^{n+1} = x_i^n + \sum_{p=0}^{m} a_{ij_0j_1\cdots j_p} \prod_{k=0}^{p} \frac{\partial f_{j_k}}{\partial x_{j_k}}\bigg|_{f_*(x_*)} \left(\overline{x_{J_k}} - x_{j_k}^n\right)$$

We arrive at the following expression for the difference between the solution $\overline{x_*}$ and the next point of iteration $x_*^{n+1}$:

$$\overline{x_\iota} - x_i^{n+1} = a_{ij_0j_1\cdots j_{m+1}} \prod_{k=0}^{m+1} \frac{\partial f_{j_k}}{\partial x_{j_k}} \left(\overline{x_{J_k}} - x_{j_k}^n\right) + \cdots \quad (22)$$

Neglecting higher order terms and substituting the expressions for the errors (18) in (22), we have:

$$\delta_i^{n+1} = a_{ij_0j_1\cdots j_{m+1}} \prod_{k=0}^{m+1} \frac{\partial f_{j_k}}{\partial x_{j_k}} \delta_{j_k}^n \quad (23)$$

Expression (23) clearly states that the order of convergence of the iteration scheme truncated to $m+2$ terms is $m+2$, since $k$ ranges from 0 to $m+1$.

The conditions for $m+2$ order of convergence are the following:

- $\det(J)|_{x_i^p} \neq 0$ for all $x_i^p$, $p = \overline{1,n}$
- $a_{ij_0j_1\cdots j_{m+1}}\big|_{x_i^p}$ is continuous for all $x_i^p$, $p = \overline{1,n}$
- $x_i^p$ is sufficiently close to the root $\bar{x}$

As for the one-dimensional case, the first condition allows for the construction of the iterative formulas, because we need to invert the Jacobian of the system to obtain the terms $a_{ij_0j_1\cdots j_m}$, see (13). The second condition ensures that the terms $a_{ij_0j_1\cdots j_q}\big|_{x_i}$ for $q = \overline{1,m}$ are computable and bounded for each iteration at every point $x_i^p$ for $p = \overline{1,n}$.

The last condition ensures $m+2$ order of convergence, because, for this to be possible, the Taylor series approximation (23) has to be valid, which involves higher order terms to be negligible.

## 6. TEST CASE

The Newton-Raphson second order iterative scheme (14), the third order Halley's method (15), and the fourth and fifth order methods, are tested in a computer algebra system called Maxima. The Newton-Raphson method is used as a comparison, while the third order extension of Halley's method [3] to systems of nonlinear equations is tested. For simplicity a two by two nonlinear system is chosen, since the actual number of variables is not important.

The system of nonlinear equations is:

$$\begin{cases} f_1(x_1, x_2) = x_1 - x_2 \\ f_2(x_1, x_2) = {x_1}^2 + x_2^2 - 2 \end{cases}$$

This system has the solutions $(1,1)$ and $(-1,-1)$.

The numerical schemes used for solving the above nonlinear system of equations are:

- 2nd order:

$$x_i^{n+1} = x_i^n - J_{k_1 i}^{-1} f_{k_1}(x_*^n)$$

- $3^{rd}$ order:

$$x_i^{n+1} = x_i^n - J_{k_1 i}^{-1} f_{k_1}(x_*^n) + \frac{1}{2} J_{k_3 k_2}^{-1} \frac{\partial J_{k_1 i}^{-1}}{\partial x_{k_2}} f_{k_1}(x_*^n) f_{k_3}(x_*^n)$$

- $4^{th}$ order:

$$\begin{aligned} x_i^{n+1} = x_i^n - J_{k_1 i}^{-1} f_{k_1}(x_*^n) &+ \frac{1}{2} J_{k_3 k_2}^{-1} \frac{\partial J_{k_1 i}^{-1}}{\partial x_{k_2}} f_{k_1}(x_*^n) f_{k_3}(x_*^n) \\ &- \frac{1}{3!} J_{k_5 k_4}^{-1} \frac{\partial}{\partial x_{k_4}} \left( J_{k_3 k_2}^{-1} \frac{\partial J_{k_1 i}^{-1}}{\partial x_{k_2}} \right) f_{k_1}(x_*^n) f_{k_3}(x_*^n) f_{k_5}(x_*^n) \end{aligned}$$

- $5^{th}$ order:

$$\begin{aligned} x_i^{n+1} = x_i^n - J_{k_1 i}^{-1} f_{k_1}(x_*^n) &+ \frac{1}{2} J_{k_3 k_2}^{-1} \frac{\partial J_{k_1 i}^{-1}}{\partial x_{k_2}} f_{k_1}(x_*^n) f_{k_3}(x_*^n) \\ &- \frac{1}{3!} J_{k_5 k_4}^{-1} \frac{\partial}{\partial x_{k_4}} \left( J_{k_3 k_2}^{-1} \frac{\partial J_{k_1 i}^{-1}}{\partial x_{k_2}} \right) f_{k_1}(x_*^n) f_{k_3}(x_*^n) f_{k_5}(x_*^n) \\ &+ \frac{1}{4!} J_{k_7 k_6}^{-1} \frac{\partial}{\partial x_{k_6}} \left( J_{k_5 k_4}^{-1} \frac{\partial}{\partial x_{k_4}} \left( J_{k_3 k_2}^{-1} \frac{\partial J_{k_1 i}^{-1}}{\partial x_{k_2}} \right) \right) f_{k_1}(x_*^n) f_{k_3}(x_*^n) f_{k_5}(x_*^n) f_{k_7}(x_*^n) \end{aligned}$$

Each of the following tables shows, for each iteration of its formula, the solution to 50 significant digits and the difference between successive solutions.

Table 1 – Newton Raphson **second order** multivariate method

| **Iteration number** | **Solution** $(x_1 = x_2)$ | **Difference between solutions** $(\delta = x_{i+1} - x_i)$ |
|---|---|---|
| **0** | 4e0 | - |
| **1** | 2.125e0 | 1.875b0 |
| **2** | 1.29779411764705882352941176470588235294117e0 | 8.272058823e-1 |
| **3** | 1.03416618063656057323779370104982502916180e0 | 2.636279370e-1 |
| **4** | 1.00056438119963058597486609415384203374824e0 | 3.360179943e-2 |
| **5** | 1.00000015917323486698635849032681600137216e0 | 5.642220263e-4 |
| **6** | 1.00000000000001266805733259473578107074834e0 | 1.591732221e-7 |
| **7** | 1.00000000000000000000000000008023983829095e0 | 1.266805733e-14 |
| **8** | 1.00000000000000000000000000000000000000000e0 | 8.023983829e-29 |
| **9** | 1.00000000000000000000000000000000000000000e0 | 3.219215824e-57 |
| **10** | 1.00000000000000000000000000000000000000000e0 | 5.181675262e-114 |
| **11** | 1.00000000000000000000000000000000000000000e0 | 1.342487926e-227 |
| **12** | 1.00000000000000000000000000000000000000000e0 | 9.011369159e-455 |
| **13** | 1.00000000000000000000000000000000000000000e0 | 4.060238706e-909 |

Table 2 – Halley's **third order** multivariate method

| **Iteration number** | **Solution** $(x_1 = x_2)$ | **Difference between solutions** $(\delta = x_{i+1} - x_i)$ |
|---|---|---|
| **0** | 4e0 | - |
| **1** | 1.685546875e0 | 2.314453125e0 |
| **2** | 1.05093669710446668578038273953086034451734e0 | 6.346101778e-1 |
| **3** | 1.00005910371154170756114074221442391204039e0 | 5.087759339e-2 |
| **4** | 1.00000000000010321825550347214780221994316e0 | 5.910371143e-5 |
| **5** | 1.00000000000000000000000000000000000000054e0 | 1.032182555e-13 |
| **6** | 1.00000000000000000000000000000000000000000e0 | 5.498440738e-40 |
| **7** | 1.00000000000000000000000000000000000000000e0 | 8.311676855e-119 |
| **8** | 1.00000000000000000000000000000000000000000e0 | 2.871018262e-355 |

Table 3 – **Fourth order** multivariate method

| Iteration number | Solution $(x_1 = x_2)$ | Difference between solutions $(\delta = x_{i+1} - x_i)$ |
|---|---|---|
| **0** | 4e0 | - |
| **1** | 1.47955322265625e0 | 2.520446777e0 |
| **2** | 1.00832805021999203253155486858343263965267e0 | 4.712251724e-1 |
| **3** | 1.00000000291805361538124559234554057497560e0 | 8.328047301e-3 |
| **4** | 1.00000000000000000000000000000000004531615e0 | 2.918053615e-9 |
| **5** | 1.00000000000000000000000000000000000000000e0 | 4.531615792e-35 |
| **6** | 1.00000000000000000000000000000000000000000e0 | 2.635677954e-138 |
| **7** | 1.00000000000000000000000000000000000000000e0 | 3.016125394e-551 |

Table 4 – **Fifth order** multivariate method

| Iteration number | Solution $(x_1 = x_2)$ | Difference between solutions $(\delta = x_{i+1} - x_i)$ |
|---|---|---|
| **0** | 4e0 | - |
| **1** | 1.358853816986083984375e0 | 2.641146183e0 |
| **2** | 1.00116069568552031665772086358127934091978e0 | 3.576931213e-1 |
| **3** | 1.00000000000000183265685289786233037850973e0 | 1.160695685e-3 |
| **4** | 1.00000000000000000000000000000000004531615e0 | 1.832656852e-15 |
| **5** | 1.00000000000000000000000000000000000000000e0 | 1.808896959e-74 |
| **6** | 1.00000000000000000000000000000000000000000e0 | 1.694639002e-369 |

From Table 1 and Table 2 it can be concluded that the number of correct digits between iterations, when the point is close to the solution, doubles, while from Table 3 and 4, the number of exact digits triples with each iteration, and so on.

This confirms that the iteration formulas built are of second, third, fourth and fifth order of convergence in the multivariable case and that formula (23) from chapter 5, which provides the order of convergence, is indeed correct.

The code used for this test case is listed below:

```
kill(all)$
numer:false$
/*============================================DEFINITIONS*/
f:[
    x[1]+2*x[2]+x[3],
    2*x[1]-x[2]-x[3],
    x[1]^2+x[2]^2+x[3]^2-3
]$
n:length(f)$
j:zeromatrix(n,n)$
for p0:1 thru n do
for p1:1 thru n do
    j[p0,p1]:diff(f[p1],x[p0])$
j:j^^-1$
id:ident(n)$
/*============================================TERMS*/
nb_term:0$

/*2nd order*/
nb_term:nb_term+1$
for i:1 thru n do
    t[nb_term,i]:sum(sum(
        j[p1,p0]*id[i,p0]*(-f[p1])
        ,p1,1,n),p0,1,n)$
```

```

/*3rd order*/
nb_term:nb_term+1$
for i:1 thru n do
    t[nb_term,i]:sum(sum(sum(sum(j[p3,p2]*diff(j[p1,p0]*id[i,p0]
    ,x[p2])*(-f[p1])*(-f[p3])/2,p3,1,n),p2,1,n),p1,1,n),p0,1,n)$

/*4th order*/
nb_term:nb_term+1$
for i:1 thru n do
    t[nb_term,i]:sum(sum(sum(sum(sum(sum(j[p5,p4]*diff(j[p3,p2]*
    diff(j[p1,p0]*id[i,p0],x[p2]),x[p4])*(-f[p1])*(-f[p3])*(-f[p
    5])/3!,p5,1,n),p4,1,n),p3,1,n),p2,1,n),p1,1,n),p0,1,n)$

/*5th order*/
nb_term:nb_term+1$
for i:1 thru n do
    t[nb_term,i]:sum(sum(sum(sum(sum(sum(sum(sum(j[p7,p6]*diff(j
    [p5,p4]*diff(j[p3,p2]*diff(j[p1,p0]*id[i,p0],x[p2]),x[p4]),x
    [p6])*(-f[p1])*(-f[p3])*(-f[p5])*(-f[p7])/4!,p7,1,n),p6,1,n)
    ,p5,1,n),p4,1,n),p3,1,n),p2,1,n),p1,1,n),p0,1,n)$

/*===========================================ITERATE*/
numer:true$
fpprec:1000$
nb_iter:10$
for i:1 thru n do(
    x0[i]:0,
    x1[i]:0,
    x2[i]:0
)$
x[1]:2.1$
x[2]:2.2$
x[3]:-1$
listarray(x);
for p0:1 thru nb_iter do(
    for i:1 thru n do(
        x_temp[i]:bfloat(x[i]),
        for p1:1 thru nb_term do
            x_temp[i]:bfloat(x_temp[i]+ev(t[p1,i]))
    ),
    for i:1 thru n do x[i]:bfloat(x_temp[i]),
    disp(""),
    display(p0),
    for i:1 thru n do(
        x0[i]:x1[i],
        x1[i]:x2[i],
        x2[i]:x[i],
        d1:bfloat(abs(x1[i]-x0[i])),
        d2:bfloat(abs(x2[i]-x1[i])),
        display(x[i]),
        disp(d2)
    )
)$
```

## 7. A BRIEF EXPLORATION OF THE STATE OF ART

In the one-dimensional case, there are extensive studies, including methods of various orders and proofs of their order convergence. Among these methods there are the second order Newton-Raphson method, the third order Halley's method, [4] and [6], the extensive work done by Schröder, [2] and [3], to encapsulate the various methods in two families of iterative schemes, Householder's method [1], and many others, [8],[5] and [7]. Regarding the multidimensional case, the most common method is the Newton-Raphson method, and other of its higher order variations. It is known that the order of convergence can be proved by considering the inverse function. To the best of our knowledge, we have not yet found papers where the idea of using multiple inverse functions, to derive and prove iterative schemes of order higher than two for systems of nonlinear equations, has been explicitly stated.

## 8. CONCLUSIONS

The detailed derivation of the iteration schemes for systems of nonlinear equations provides a simple way of quickly constructing such schemes when a higher order of convergence is needed. The idea of using multiple inverse functions explains some of the general structure of iterative numerical schemes, and also lends more insight to the user over their properties. The generalised proof, presented in chapters three and five, offers a quick way to check the order of convergence of any iterative numerical scheme used for solving nonlinear systems of equations.

The original contributions of the author are:

- A quick methodology to construct iterative numerical schemes of any order for solving nonlinear systems of equations.
- Formulas for calculating the order of convergence of the iterative numerical schemes constructed using the methodology presented in this paper.
- Insights regarding the form and application of the iterative schemes.

## REFERENCES

[1] A. S. Householder, *The Numerical Treatment of a Single Nonlinear Equation*, New York: McGraw-Hill, 1970.

[2] E. Schröder, Über unendlich viele Algorithmen zur Auflösung der Gleichungen, *Math. Ann.* 2, 317-365, 1870.

[3] G. W. Stewart, *On Infinitely Many Algorithms for Solving Equations*, English translation of Schröder's original paper, College Park, MD: University of Maryland, Institute for Advanced Computer Studies, Department of Computer Science, 1993.

[4] T. R Scavo and J. B. Thoo, On the Geometry of Halley's Method, *Amer. Math.* Monthly **102**, 417-426, 1995.

[5] J. R. Sharma. A Family of Third-Order Methods to Solve Nonlinear Equations by Quadratic Curves Approximation, *Applied Mathematics and Computation*, **184**: 210–215, 2007.

[6] M. A. Hernandez and M. A. Salanova, A Family of Chebyshev–Halley Type Methods, *International Journal of Computer Mathematics*, **47**: 59–63, 1993.

[7] J. M. Ortega, and W. G. Rheinboldt, *Iterative Solution of Nonlinear Equations in Several Variables*, Academic Press, New York, 1970.

[8] B. T. Polyak, Newton's Method and Its Use in Optimization. *European Journal of Operational Research*, **127**(3): 1086–1096, 2007.